\documentclass [english,a4paper,10pt]{article}
\usepackage {babel}
\usepackage {array,theorem}
\usepackage {amstex,amscd}

\newtheorem{theorem}{Theorem}[section]
\newtheorem{lemma}[theorem]{Lemma}
\newtheorem{proposition}[theorem]{Proposition}
\newtheorem{cor}[theorem]{Corollary}
\theorembodyfont{\rm}
\newtheorem{remark}[theorem]{Remark}
\newtheorem{remarks}[theorem]{Remarks}
\newtheorem{definition}{Definition} 
\newtheorem{uremarks}{Remarks}      
        
\newtheorem{claim}{Claim}        

\newcommand{\PP}{{\Bbb{P}}}

\newcommand{\calO}{{\cal O}}

\newcommand{\calI}{{\cal I}}

\newcommand{\on}[1]{\operatorname{#1}}
\renewcommand{\phi}{\varphi}
\newenvironment{Proof}{\begin{ProofwCaption}{Proof}}{\end{ProofwCaption}}
\newenvironment{Proof*}[1]{\begin{ProofwCaption}{{#1}}}{\end{ProofwCaption}}
\newenvironment{ProofwCaption}[1]%
  {\addvspace\theorempreskipamount \noindent{\it #1.}\rm}%
  {\qed \par \addvspace\theorempostskipamount}
\newcommand{\qedsymbol}{\mbox{$\Box$}}
\newcommand{\qed}{\quad\qedsymbol}
\begin{document}
\title{Rational surfaces in $\PP^4$ containing a plane curve}
\author{F.~Catanese and K.~Hulek}
\markboth{Catanese/Hulek}{Rational surfaces in $\PP^4$ containing a plane
curve}
\date{}
\maketitle

\begin{abstract}
The families of smooth rational surfaces in $\PP^4$ have been
classified in degree $\le 10$. All known rational surfaces in $\PP^4$ can be
represented as blow-ups of the plane $\PP^2$. The fine classification of these
surfaces consists of giving explicit open and closed conditions which determine
the configurations of points corresponding to all surfaces in a given family.
Using a restriction argument originally due independently to Alexander and
Bauer we achieve the fine classification in two cases, namely non-special
rational surfaces of degree 9 and special rational surfaces of degree 8. The
first case completes the fine classification of all non-special rational
surfaces. In the second case we obtain a description of the moduli space as the
quotient of a rational variety by the symmetric group $S_5$. We also discuss in
how far this method can be used to study other rational surfaces in $\PP^4$.
\end{abstract}

\section{Introduction}\label{sectionI}

The families of smooth rational surfaces in $\PP^4$ have been
classified in degree
$\leq10$ ( \cite{A1}, \cite{I1}, \cite{I2}, \cite{O1}, \cite{O2},
\cite{R1}, \cite{R2}, \cite{PR}). In this thesis Popescu
\cite{P} constructed further examples of rational surfaces in degree $11$. The
existence of these surfaces has been proved in various ways, using linear
systems, vector bundles and sheaves or liaison arguments. All known rational
surfaces can be represented as a blowing-up of $\PP^2$. Although it would
seem the most natural approach to prove directly that a given linear system is
very ample, this turns out to be a very subtle problem in some cases, in
particular when the surface $S$ in $\PP^4$ is special
(i.e.~$h^1(\calO_S(H))\neq0$). On the other hand, being able to handle the
linear system often means that one knows the geometry of the surface very well.

The starting point of our paper is the observation that every known rational
surface in $\PP^4$ contains a plane curve $C$. Using the hyperplanes through
$C$ one can construct a residual linear system $|D|$. I.e., we can write
$H\equiv C+D$ with $\dim |D|\geq1$. This situation was studied in particular
by Alexander \cite{A1}, \cite{A2} and Bauer \cite{B}: if $|H|$ restricts to
complete linear systems on $C$ and $D'$ where $D'$ varies in a $1$-dimensional
linear subsystem of $|D|$, then $H$ is very ample on $S$ if and only if it is
very ample on $C$ and the curves $D'$ (cf.~Theorem~(\ref{theo21})). In this
way one can reduce the question of very ampleness of $H$ to the study of
linear systems on curves. In \cite{CFHR} the following curve embedding theorem
was proved which we shall state here only for the (special) case of curves
contained in a smooth surface.

\begin{theorem}\label{theo11}
A divisor $H$ is very ample on $C$ if for every subcurve $Y$ of $C$ of
arithmetic genus $p(Y)$

\noindent $\on{(i)}$ $H.Y\geq 2p(Y)+1$ or

\noindent $\on{(ii)}$ $H.Y\geq 2p(Y)$ and there is no $2$-cycle $\xi$ of $Y$
such that $I_\xi\calO_Y\cong \omega_Y(-H)$.

More generally

\noindent $\on{(iii)}$ If $\xi$ is an $r$-cycle of $C$, then
$H^0(C,\calO_C(H))$ surjects onto $H^0(\calO_C(H)\otimes\calO_\xi)$ unless
there is a subcurve $Y$ of $C$ and a morphism $\phi:I_\xi\calO_Y\to
\omega_Y(-H)$ which is "good" (i.e.~$\phi$ is injective with a cokernel of
finite length) and which is not induced by a section of $H^0(Y,\omega_Y(-H))$.
\end{theorem}

The method described above was used in \cite{CF} to characterize exactly all
configurations of points in $\PP^2$ which define non-special rational surfaces
of degree $\leq8$. In these cases $H.D\geq 2p(D)+1$. This left the case open
of one non-special surface, namely the unique non-special surface of degree
$9$. In this case one has a decomposition $H\equiv C+D$ where $C$ is a plane
cubic, and $|D|$ is a pencil of curves of genus $p(D)=3$ and $H.D=6$.
Section~\ref{sectionII} is devoted to this surface. In Theorem~(\ref{theo22})
we classify all configurations of points in the plane which lead to
non-special surfaces of degree $9$ in $\PP^4$. This completes the fine
classification of non-special surfaces.

In section~\ref{sectionIII} we show that this method can also be applied to
study special surfaces. We treat the (unique) special surface of degree $8$.
In this case there exists a decomposition $H\equiv C+D$ where $C$ is a conic
and $|D|$ is a pencil of curves of genus $4$ with $H.D=6$. It turns out that
for the general element $D'$ of $|D|$ (but not necessarily for all elements)
$H$ is the canonical divisor on $D'$. In Theorem~(\ref{IIItheo14}) we give a
characterization of these configurations of points which define smooth special
surfaces of degree $8$ in $\PP^4$. We then use this result to give an
existence proof (in fact we construct the general element in the family) of
these surfaces using only the linear system $|H|$
(Theorem~(\ref{IIItheo17})), and in particular to describe the moduli
space of the above surfaces modulo projective
equivalence (Theorem~(\ref{IIItheo20})).

Finally in section~\ref{sectionIV} we discuss some posibilities how this
method can be used to study other rational surfaces in $\PP^4$,
suggesting some explicit decompositions $H\equiv C+D$ of the hyperplane
class as the sum of divisors.

\medskip

\noindent {\bf Acknowledgements.} The authors were partially supported by the
DFG-Schwerpunktprogramm "Komplexe Mannigfaltigkeiten" under contract number
Hu337/4-3, the EU HCM project AGE (Algebraic Geometry in Europe) contract
number ERBCHRXCT 940557 and MURST $40\%$. The second author is also
grateful to the Tata-Institute of Fundamental Research for their
hospitality. The final version was written while the first author was
"{\em Professore distaccato}" at the {\em Accademia dei Lincei}.

\section{The non-special rational surface of degree $9$}\label{sectionII}

In this section we want to give an application of Theorem (\ref{theo11})
to non-special rational surfaces. These surfaces have been classified by
Alexander \cite{A1}. Catanese and Franciosi treated all
non-special rational surfaces of degree $\leq8$ by studying suitable
decompositions $H=C+D$ of the embedding linear systems. The crucial
observation here is the following result, originally due to J.~Alexander and
I.~Bauer \cite{B}.

\begin{theorem}[Alexander-Bauer]\label{theo21}
Let $X$ be a smooth projective variety and let $C,D$ be effective divisors with
$\dim|D|\geq1$. Let $H$ be the divisor $H\equiv C+D$. If $\big|H\big||_C$ is
very ample and for all $D'$ in a $1$-dimensional subsystem of $|D|$,
$\big|H\big||_{D'}$ is very ample, then $|H|$ is very ample on $X$.
\end{theorem}

By Alexander's list there is only one non-special rational surface of
degree bigger than $8$. This surface is a $\PP^2$ blown up in 10 points
$x_1,\ldots,x_{10}$ embedded by the linear system $|H|=|13L-4
\sum_{i=1}^{10}x_i|$. Alexander showed that for general position of the
points $x_i$ the linear system $|H|$ embeds
$S={\tilde{\PP}}^2(x_1,\ldots,x_{10})$ into $\PP^4$. Clearly the degree
of $S$ is 9. Here we show that using Theorem (\ref{theo11}) one can also
apply the decomposition method to this surface. In fact we obtain
necessary and sufficient conditions for the position of the points
$x_i$ for $|H|$ to be very ample. Our result is the following

\begin{theorem}\label{theo22} The linear system $|H|= |13L-4\sum x_i|$
embeds the surface
$S ={\tilde \PP}^2(x_1,\ldots, x_{10})$ into
$\PP^4$ if and only if

\noindent $\on{(0)}$ no $x_i$ is infinitely near,

\noindent $\on{(1)}$ $|L-\sum\limits_{i\in\triangle}x_i|=\emptyset$ for
$|\triangle|\geq4$,

\noindent $\on{(2)}$ $|2L-\sum\limits_{i\in\triangle'}x_i|=\emptyset$ for
$|\triangle'|\geq7$,

\noindent $\on{(3)}$ $|3L-\sum\limits_ix_i|=\emptyset$,

\noindent $\on{(3)}_{ij}'$ $|3L-\sum\limits_{k\neq
i,j}x_k-2x_i|=\emptyset$ for all pairs $(i,j)$,

\noindent $\on{(4)}_{ijk}$ $|4L-2x_i-2x_j-2x_k- \sum\limits_{l\neq
i,j,k}x_l |=\emptyset$ for all triples $(i,j,k)$,

\noindent $\on{(6)}_i$ $|6L-x_i-2\sum\limits_{j\neq i}x_j|=\emptyset$,

\noindent $\on{(10)}_1$  If $D=10L-4x_1-3\sum\limits_{i\geq2}x_i$, then
$\dim|D|=1$.
\end{theorem}

\begin{uremarks} (i) Clearly conditions (0) to (6) are open conditions.
The expected dimension of $|D|$ is 1, hence this condition is also open.

\noindent (ii) The last condition is asymmetrical. If $|H|$ is very ample
condition $(10)_i$ is necessarily fulfilled for all $i$. On the other
hand, our theorem shows that in order to prove very ampleness for $|H|$ it
suffices to check only one of the conditions $(10)_i$.
\end{uremarks}

\begin{Proof} We shall first show that the conditions stated are
necessary. Clearly (0) follows since $H.(x_i-x_j)=0$. Similarly the
ampleness of $H$ immediately implies conditions (1) to (4). Assume the
linear system
$|6L-x_i-2\sum\limits_{j\neq i}x_j|$ contains some element $A$. Then
$H.A=2$, and $p(A)=1$ which contradicts very ampleness of $H$. For (10) we
consider $C\equiv H-D\equiv 3L-\sum\limits_{i\geq2}x_i$. Clearly $|C|$ is
non empty. For $C'\in|C|$ we consider the exact sequence
\setcounter{equation}{10}
\begin{equation}\label{gl11}
0\longrightarrow\calO_S(D)\longrightarrow\calO_S(H)\longrightarrow
\calO_{C'}(H)\longrightarrow 0.
\end{equation} If $h^0(\calO_S(D))\geq3$, then either
$h^0(\calO_S(H)\geq6$ and $|H|$ does not embed $S$ into $\PP^4$ or $|H|$
maps $C'$ to a line. But since $p(C)=1$ this means that $|H|$ cannot be
very ample.

Now assume that conditions (0) to $(10)_1$ hold. We shall first show
\begin{align*} h^1(\calO_S(D))&=0\tag{I} \\ h^1(\calO_S(C))&=0\tag{II}\\
h^0(\calO_S(H))&=5 \tag{III}
\end{align*} Ad (I): By condition $(10)_1$ we have $h^0(\calO_S(D))=2$.
Clearly
$h^2(\calO_S(D))=h^0(\calO_S(K-D))=0$. Hence the claim follows from
Riemann-Roch, since $\chi(\calO_S(D))=2$.

\noindent Ad (II): We consider $-K\equiv 3L-\sum\limits_ix_i\equiv
C-x_1$. By condition (3) $h^0(\calO_S(-K))=0$. Clearly also
$h^2(\calO_S(-K))=h^0(\calO_S(K))=0$. Hence by Riemann-Roch
$h^1(\calO_S(-K))=-\chi (\calO_S(-K))=0$. Now consider the exact sequence
\begin{equation}\label{gl12}
0\longrightarrow\calO_S(-K)\longrightarrow\calO_S(C)\longrightarrow
\calO_S(C)|_{x_1}=\calO_{x_1}\longrightarrow 0.
\end{equation} This shows $h^1(\calO_S(C))=0$. Note that this also
implies (by Riemann-Roch) that $h^0(\calO_S(C))=1$, i.e.~the curve $C'$
is uniquely determined.

\noindent Ad (III): In view of (I) and sequence (\ref{gl11}) it suffices
to show that $h^0(\calO_{C'}(H))=3$. By Riemann-Roch on $C'$ this is
equivalent to
$h^1(\calO_{C'}(H))=0$. Since $K_{C'}$ is trivial this in term is
equivalent to
$h^0(\calO_{C'}(-H))=0$. By condition $(3)'$ the curve $C'$ contains no
exceptional divisor. As a plane curve $C'$ can be irreducible or it can
decompose into a conic and a line or three lines. In view of conditions
(1) and (2), however, $C'$ cannot have multiple components and, moreover
$H$ has positive degree on every component. This proves
$h^0(\calO_{C'}(-H))=0$ and hence the claim.

This shows that $|H|$ maps $S$ to $\PP^4$ and that, moreover, $|H|$
restricts to complete linear systems on $C'$ and all curves $D'\in|D|$.
We shall now show
\begin{align*}
\begin{split} &\ \ \ \,\text{ For every subcurve }A\leq C'\text{ we have }
H.A\geq2p(A)+1
\end{split}\tag{IV}
\end{align*}
\begin{align*}
\begin{split} &\text{For every proper subcurve }B'\subset D'\text{ of an
element}\\ &D'\in |D|\text{ we have }H.B'\geq2p(B')+1
\end{split}\tag{V(i)}
\end{align*}
\begin{align*}
\begin{split} &H\text{ does not restrict to a "$(2+K)$"-divisor on }D',\\
&\text{i.e. }\calO_{D'}(H-K_{D'})\text{ does not have a good section}\\
&\text{defining a degree } 2\text{-cycle}.
\end{split}\tag{V(ii)}
\end{align*}
It then follows from (IV) and \cite[Theorem 3.1]{CF} that
$|H|$ is very ample on $C'$. Because of (V(i)) and (V(ii)) it follows
from Theorem (\ref{theo11}) that $|H|$ is very ample on every element
$D'$ of $|D|$. It then follows from Theorem (\ref{theo21}) that $|H|$ is
very ample.

\noindent Ad (V(ii)): Let $H_{D'}$ be the restriction of $H$ to $D'$, and
denote the canonical bundle of $D'$ by $K_{D'}$. It suffices to show that
$h^0(\calO_{D'}(H_{D'}-K_{D'}))=0$. Now
\begin{align*} H_{D'}-K_{D'}&= (H-K-D)|_{D'}\\ &= (C-K)|_{D'}\\ &=
(2C-x_1)|_{D'}.
\end{align*} There is an exact sequence
\begin{equation}\label{gl13}
0\longrightarrow\calO_S(2C-x_1-D)\longrightarrow\calO_S(2C-x_1)
\longrightarrow
\calO_{D'}(H_{D'}-K_{D'})\longrightarrow 0.
\end{equation} Since
$$ 2C-x_1\equiv 6L-x_1-2\sum_{i=2}^{10}x_i
$$ it follows from condition $(6)_1$ that $h^0(\calO_S(2C-x_1))=0$.
Clearly
$h^0(\calO_S(2C-x_1-D))=0$. Now
$$ 2C-x_1-D\equiv -4L+3x_1+\sum_{i=2}^{10}x_i
$$ resp.
$$ K-(2C-x_1-D)\equiv L-2x_1.
$$ Hence $h^2(\calO_S(2C-x_1-D))=h^0(\calO_S(K-(2C-x_1-D))=0$. Since
moreover
$\chi(\calO_S(2C-x_1-D))=0$ it follows that $h^1(\calO_S(2C-x_1-D))=0$.
The assertion follows now from sequence (\ref{gl13}).

\noindent Ad (IV) and (V(i)): We have to show that for all curves $A$
with $A\leq C'$, resp.~$A<D'$, $D'\in |D|$ the following holds
\begin{equation}\label{gl14} H.A\geq2p(A)+1.
\end{equation} We first notice that it is enough to prove (\ref{gl14})
for divisors $A$ with
$p(A)\geq0$. Assume in fact we know this and that $p(A)<0$. Then $A$ is
necessarily reducible. For every irreducible component $A'$ of $A$ we
have $p(A')\geq0$ and hence $H.A'>0$. This shows $H.A>0$ and hence
(\ref{gl14}). Clearly (\ref{gl14}) also holds for the lines $x_i$. Hence
we can assume that $A$ is of the form
\begin{equation}\label{gl15} A\equiv aL-\sum_ib_ix_i\quad\text{with
}1\leq a\leq 10.
\end{equation} Note that
\begin{eqnarray} 2p(A)&=&a(a-3)-\sum_ib_i(b_i-1)+2\label{gl16}\\
H.A&=&13a-4\sum_ib_i.\label{gl17}
\end{eqnarray} We proceed in several steps

\medskip

\noindent {\bf Claim 1 } Let $A$ be as in (\ref{gl15}) with $1\leq a\leq
3$. Assume that $p(A)\geq0$. Then (\ref{gl14}) is fulfilled.

\medskip

\noindent {\em Proof of Claim 1 } After possibly relabelling the $x_i$ we
can assume that $b_1\geq b_2\geq\ldots\geq b_{10}$. If $a=1$ or 2 then
$b_1\leq1$ and $b_2\geq0$. Moreover $p(A)=0$. If $H.A\leq2p(A)$ we get
immediately a contradiction to conditions (1) or (2). If $a=3$ then we
have two cases. Either
$b_1\leq1$, $b_2\geq0$ as above and $p(A)=1$. Then $H.A\leq2p(A)$ violates
condition (3). Or $b_1=2$ or $b_{10}=-1$ and the other $b_i$ are 0 or 1.
Then
$H.A\leq2p(A)$ is only possible for $b_1=2$, but this would violate
condition
$(3)'$.

\medskip

\noindent {\bf Claim 2 } $H$ is ample on $C$ and $D$, i.e.~for every
irreducible component $A$ of $C'$, resp.~$D'$, $D'\in|D|$ we have $H.A>0$.

\medskip

\noindent {\em Proof of Claim 2 } Assume the claim is false. Let $A$ be an
irreducible component with $H.A\leq0$. Since $A$ is irreducible,
$p(A)\geq0$. By (\ref{gl16}), (\ref{gl17}) this leads to the two
inequalities
\begin{eqnarray} 13a&\leq&4\sum b_i\label{gl18}\\
\sum b_i(b_i-1)&\leq&a(a-3)+2.\label{gl19}
\end{eqnarray} Multiplying (\ref{gl19}) by $13^2$ and using (\ref{gl18})
we obtain
\begin{equation}\label{gl20} 169\left(\sum b_i^2-\sum b_i\right)\leq
16\left(\sum b_i\right)^2- 156\sum b_i+338.
\end{equation} Now
\begin{equation}\label{gl21}
\left(\sum b_i\right)^2=10\sum b_i^2-\sum_{i<j}(b_i-b_j)^2
\end{equation} and using this (\ref{gl20}) becomes
\begin{equation}\label{gl22}
\sum_i(9b_i^2-13b_i)+16\sum_{i<j}(b_i-b_j)^2\leq 338.
\end{equation} The function $f(b)=9b^2-13b$ for integers $b$ is non
positive only for $b=0$ or 1. It is minimal for $b=1$. Since $f(1)=-4$ we
derive from (\ref{gl22})
\begin{equation}\label{gl23} 16\sum_{i<j}(b_i-b_j)^2\leq 378
\end{equation} resp.
\begin{equation}\label{gl24}
\sum_{i<j}|b_i-b_j|^2\leq 23.
\end{equation} At this point it is useful to introduce the following
integer valued function
$$
\delta=\delta(A)=\max_{i<j}|b_i-b_j|.
$$ We have to distinguish several cases:

\noindent $\delta\geq3$: Assume there is a pair $(i,j)$ with
$|b_i-b_j|\geq3$. Then for all $k\neq i,j$:
$$ |b_i-b_k|^2+|b_j-b_k|^2\geq5.
$$ Hence
$$
\sum_{i<j}|b_i-b_j|^2\geq9+5\cdot 8=49
$$ contradicting (\ref{gl24}).

\noindent $\delta=2$: After possibly relabelling the $x_i$ we can assume
that
$b_2=b_1+2$ and $b_1\leq b_k\leq b_2$ for $k\geq3$. Then
$$ |b_k-b_1|^2+|b_k-b_2|^2=\left\{\begin{array}{cl} 2&\text{if
}b_k=b_1+1\\ 4&\text{if }b_k=b_1\text{ or }b_k=b_2. \end{array}\right.
$$ Let $t$ be the number of $b_k$ which are either equal to $b_1$ or
$b_2$. Then
\begin{eqnarray*}
\sum_{i<j}|b_i-b_j|^2&\geq&4+4t+2(8-t)+t(8-t)\\ &=&20+t(10-t).
\end{eqnarray*} It follows from (\ref{gl24}) that $t=0$. But then
(\ref{gl22}) gives
$$
\sum(9b_i^2-13b_i)\leq18.
$$ Looking at the values of $f(b)=9b^2-13b$ one sees immediately that
this is only possible for $b_1=-1$ or $b_1=0$. In the first case it
follows from (\ref{gl18}) that $a<0$ which is absurd. In the second case
we obtain $a\leq 3$ and hence we are done by Claim 1.

\noindent $\delta\leq1$: Here we can assume
$$ b_1=\ldots=b_k=m,\quad b_{k+1}=\ldots=b_{10}=m+1.
$$ Since $f(b)\geq42$ for $b\geq3$ it follows immediately from
(\ref{gl22}) that
$m\leq2$. If $m\leq0$ then (\ref{gl18}) gives $a\leq3$ and we are done by
Claim 1. It remains to consider the subcases $m=1$ or 2.

\noindent $m=2$: Since $f(2)=10$ and $f(3)=42$ formula (\ref{gl22})
implies
$$ 10k+42(10-k)+16k(10-k)\leq338.
$$ One checks easily that this is only possible for $k=9$ or 10. In this
case (\ref{gl18}) gives $a\leq6$. If $k=9$ then (\ref{gl18}) gives
$22\leq a(a-3)$, i.e.~$a\geq7$, a contradiction. If $k=10$, then
(\ref{gl18}) implies
$18\leq a(a-3)$. This is only possible for $a=6$. But now the existence
of $A$ would contradict condition (6).

\noindent $m=1$: Since $f(1)=-4$ and $f(2)=10$ formula (\ref{gl22}) reads
$$ -4k+10(10-k)+16k(10-k)\leq338
$$ or equivalently
$$ k(73-8k)\leq119.
$$ It is straightforward to check that this implies $k\leq2$ or $k\geq7$.
If
$k\leq2$ then $\sum b_i(b_i-1)\geq16$ and (\ref{gl19}) shows that
$a\geq6$. On the other hand $\sum b_i\leq19$ and this contradicts
(\ref{gl18}). Now assume
$k\geq7$. Then $\sum b_i\leq13$. It follows from (\ref{gl18}) that either
$a\leq3$ -- and this case is dealt with by Claim 1 -- or $a=4$ and $\sum
b_i=13$. Then $k=7$ and the existence of $A$ contradicts condition (4).

\medskip

\noindent {\em End of proof } It follows immediately from Claim 1 that
(\ref{gl14}) holds for subcurves $A\leq C'$. It remains to consider
subcurves
$A< D'$, $D'\in|D|$. Since $H$ is ample on $D$ we have $H.A>0$, hence it
suffices to consider curves with $p(A)\geq1$. Also by ampleness of $H$ on
$D$ it follows that
\begin{equation}\label{gl25} 1\leq H.A\leq5
\end{equation} since $H.D=6$. Also note that, as an immediate consequence
of (\ref{gl17}):
\begin{equation}\label{gl26} a\equiv H.A\on{mod}4.
\end{equation} Finally we remark the following

\medskip

\noindent {\bf Observation:} If $A<D$ is not one of the exceptional lines
$x_i$, then $H.A\leq4$ implies $b_i\geq0$ for all $i$. Otherwise at most
one
$b_i=-1$ and all other $b_i\geq0$.

This follows from the ampleness of $H$ on $x_i$, since $H.x_i=4$.

{}From now on we set
\begin{equation}\label{gl27} B:=D-A.
\end{equation} By adjunction
\begin{equation}\label{gl28} p(A)+p(B)=p(D)+1-A.B=4-A.B.
\end{equation} We write
$$ B\equiv bL-\sum c_ix_i.
$$ We shall now proceed by discussing the possible values of the
coefficient $a$ of $A$ in decreasing order.

\noindent $a=10$: Then $B=\sum c_ix_i$, $c_i\geq0$ and since $H.B\leq5$
we must have $B=x_i$. Then $A.B=4$ or 5 and $p(A)\leq0$ by (\ref{gl28}).

\noindent $a=9$: By (\ref{gl25}), (\ref{gl26}) we have to consider two
cases
\begin{align*}
H.A=5&, H.B=1 \tag{$\alpha$}\\ H.A=1&, H.B=5.\tag{$\beta$}
\end{align*}
Using our above observation for $B$ in case $(\alpha)$ we find that
$$
B\equiv L-x_i-x_j-x_k.
$$
But now $A.B\geq2$ and hence $p(A)\leq1$. Hence $H.A=5\geq2p(A)+1$.

Using condition (1) we have to consider the following cases for $(\beta)$:
\begin{eqnarray*} B&\equiv&L-x_i-x_j\\ B&\equiv&L-x_i-x_j-x_k+x_l.
\end{eqnarray*} In the first case $A.B\geq4$ and $p(B)=0$, hence
$p(A)\leq0$. In the second case $A.B\geq5$ and $p(B)=-1$, hence again
$p(A)\leq0$.

\noindent $a=8$: Here the only possibility is
$$
H.A=4,\quad H.B=2.
$$
Using our observation for $B$ we find that
$$
B\equiv 2L-x_{i_1}-\ldots-x_{i_6}.
$$
Either the $x_{i_j}$ are all different or we have 1 double point (and
$B$ is a pair of lines) or 3 double points (and $B$ is a double line).
Then $A.B\geq3$ (resp.~4, resp.~8) and $p(B)=0$ (resp.~$-1$, resp.~$-3$).
In either case
$p(A)\leq1$ and hence $H.A\geq2p(A)+1$.

\noindent $a=7$: In this case
$$ H.A=H.B=3.
$$ All coefficients $b_i\geq0$. It is enough to consider divisors $A$ with
$p(A)\geq2$. Together with $H.A=3$ this leads to the following conditions
on the $b_i$:
$$
\sum b_i=22,\quad \sum b_i(b_i-1)\leq26.
$$ Let $\beta_i=\max(0,b_i-1)$. Then these conditions become
$$
\sum\beta_i\geq12,\quad\sum(\beta_i+\beta_i^2)\leq26
$$ and it is easy to check that no solutions exist.

\noindent $a=6$: We now have to consider
$$ H.A=2,\quad H.B=4.
$$ We have to consider divisors $A$ with $p(A)\geq1$. Arguing as in the
case $a=7$ this leads to
$$
\sum b_i=19,\quad \sum b_i(b_i-1)\leq18
$$ resp.
$$
\sum\beta_i\geq9,\quad\sum(\beta_i+\beta_i^2)\leq18.
$$ The only solution is $b_j=1$ for one $b_j$ and $b_i=2$ for $j\neq i$.
But then
$A\in|6L-x_j-2\sum\limits_{i\neq j}x_i|$ contradicting condition (6).

\noindent $a=5$: Then we have two possible cases
\begin{align*} H.A=5&, H.B=1 \tag{$\alpha$}\\ H.A=1&, H.B=5.\tag{$\beta$}
\end{align*} We shall treat $(\alpha)$ first. Then by the ampleness of
$H$ the curve $B$ must be irreducible. Set
$$ B=5L-\sum c_ix_i,\quad c_i\geq0.
$$ Then $H.B=1$ and irreducibility of $B$ gives:
$$
\sum c_i=16,\quad \sum c_i(c_i-1)\leq12.
$$ One easily checks that this is only possible if 6 of the $c_i$ are 2,
and the others are 1. Hence
$$ B\in|5L-2\sum_{i\in\triangle}x_i-\sum_{i\not\in\triangle}x_i|,\quad
|\triangle|=6.
$$ Then $p(B)=0$. Moreover $A.B\geq3$, hence $p(A)\leq1$ and hence
$H.A\geq2p(A)+1$.

In case $(\beta)$ we apply the above argument to $A$ and find $p(A)=0$,
i.e.~again $H.A\geq2p(A)+1$.

\noindent $a=4:$ Then $H.A=4$ and $H.B=2$. We are done if $p(A)\le 1$,
and otherwise $H.A\ge 52-44=8$, a contradiction.

\noindent $1\leq a\leq 3$: This follows immediately from Claim 1.

\noindent $a=0$: The only possibility is $A=x_i$ when nothing is to show.

This finishes the proof of the theorem.
\end{Proof}

\section{The special rational surface of degree $8$ in $\PP^4$}
\label{sectionIII}

In this section we want to show how the decomposition method can be employed
to obtain very precise geometric information also about special surfaces. We
consider the rational surface in $\PP^4$ of degree $8$, sectional genus
$\pi=6$ and speciality $h=h^1(\calO_S(1))=1$. This surface was first
constructed by Okonek \cite{O2} using reflexive sheaves. In geometric
terms it is $\PP^2$ blown-up in $16$ points embedded by a linear system of
the form
$$
|H|=|6L-2\sum_{i=1}^4x_i-\sum_{k=5}^{16}y_k|.
$$
Our aim is to study the precise open and closed conditions which the points
$x_i,y_k$ must fulfill for $|H|$ to be very ample. If $|H|$ is very ample, the
exceptional lines $x_i$ are mapped to conics. Their residual intersection with
the hyperplanes gives a {\em pencil} $|D_i|$. Hence we immediately obtain the
(closed) necessary condition
\begin{gather}
|D_i|\equiv |6L-3x_i-2\sum_{j\neq i}x_j-\sum_{k=5}^{16}y_k|\text{ is a pencil}
\tag{$D_i$}
\end{gather}
By Riemann-Roch this is equivalent to $h^1(\calO_S(D_i))=1$. We first want to
study the linear system $|H|$ on the elements of the pencil $|D_i|$. Note that
$$
p(D_i)=4,\ H.D_i=6.
$$
If $D=A+B$ is a decomposition of some element $D\in |D_i|$, then
\begin{gather}
p(A)+p(B)+A.B=5\label{IIIgl1}\\
A.H+B.H=6.\label{IIIgl2}
\end{gather}
The first equality can be proved by adjunction, the second is obvious.

\begin{lemma}\label{IIIlemma1}
Assume $|H|$ is very ample. Then for every proper subcurve $Y$ of an element
$D\in |D_i|$, $h^1(\calO_Y(H))\leq1$ and $p(Y)\leq3$.
\end{lemma}
\begin{Proof}
Riemann-Roch on $Y$ gives
\begin{gather}
h^0(\calO_Y(H))=h^1(\calO_Y(H))+H.Y+1-p(Y).\label{IIIgl3}
\end{gather}
Consider the sequence
\begin{gather}
0\longrightarrow\calO_S(H-Y)\longrightarrow\calO_S(H)
\overset{\alpha}{\longrightarrow} \calO_Y(H)
\longrightarrow 0.
\label{IIIgl4}
\end{gather}
Since $h^2(\calO_S(H-Y))=h^0(\calO_S(K-(H-Y)))=0$ and $h^1(\calO_S(H))=1$ we
have $h^1(\calO_Y(H))\leq1$. We now consider the rank of the restriction map
$H^0(\alpha)$. Since $Y$ is a curve contained in a hyperplane section
$2\leq\on{rank}(\alpha)\leq4$. If $\on{rank}\alpha=2$, then $Y$ is a
line, hence $p(Y)=0$. Next assume $\on{rank}(\alpha)=3$. In this case $Y$ is a
plane curve of degree $d=Y.H$. Since $Y$ is a proper subcurve of $D$ which is
not a line $2\leq d\leq 5$. Then $h^1(\calO_Y(H))=h^0(\calO_{\PP^2}(d-4))$.
Since $h^1(\calO_Y(H))\leq1$ this shows in fact $d\leq4$. But then
$p(Y)\leq3$. Finally assume that $\on{rank}(\alpha)=4$, i.e.~$Y$ is a space
curve. By (\ref{IIIgl3})
$$
p(Y)=h^1(\calO_Y(H))-h^0(\calO_Y(H))+H.Y+1\leq3
$$
since $H.Y\leq5$.
\end{Proof}

\begin{remark}\label{IIIrem2}
Note that the above proof also shows the following: If $Y$ is a proper
subcurve of $D$ with $p(Y)=3$, then $Y$ is a plane quartic with $H_Y=K_Y$ or
$Y$ has degree $5$.
\end{remark}

Before proceeding we note the following result from \cite{CF} which we shall
use frequently in the sequel.

\begin{proposition}\label{IIIprop3}
Let $Y$ be a curve contained in a smooth surface with $p(Y)\leq2$. If $H$
is very ample on
$S$, then
$H.Y\geq 2p(Y)+1$.
\end{proposition}
\begin{Proof}
\cite[Prop.~5.2]{CF}
\end{Proof}

\begin{proposition}\label{IIIprop4}
If $|H|$ is very ample, then every element $D\in |D_i|$ is $2$-connected.
Moreover, either

\noindent $\on{(i)}$ $D$ is $3$-connected or

\noindent $\on{(ii)}$ Every decomposition of $D$ which contradicts
$3$-connectedness is either of the form $D=A+B$ with $H.B=4$, $H_B=K_B$ or of
the form $D=A+B$ with $H.B=5$. In the latter case $B=B'+B''$ with $H.B'=4$,
$H_{B'}=K_{B'}$.
\end{proposition}
\begin{Proof}
Let $D=A+B$. We first consider the case $p(A),p(B)>0$. Since $|H|$ is very
ample, it follows that $H.A\geq3$, $H.B\geq3$. But then $H.A=H.B=3$ and hence
$p(A)=p(B)=1$. By (\ref{IIIgl1}) this shows $A.B=3$.

Now assume $p(A)\leq0$. Since $p(B)\leq3$ by Lemma~(\ref{IIIlemma1}) it
follows from (\ref{IIIgl1}) that $A.B\geq2$. The only case where $A.B=2$ is
possible is $p(A)=0$, $p(B)=3$. In this case $H.B\geq4$ since Riemann-Roch for
$B$ gives
$$
h^0(\calO_B(H))=h^1(\calO_B(H))+H.B-2
$$
and we know that $h^0(\calO_B(H))\geq3$. We first treat the case $H.B=4$. Then
$h^1(\calO_B(H))=1$ and $h^0(\calO_B(H))=3$. In this case $B$ is a plane
quartic and $H_B=K_B$. Now assume $H.B=5$. If $h^1(\calO_B(H))=0$ then $B$ is
a plane quintic. But in this case $p(B)=6$, a contradiction. It remains to
consider the case $h^1(\calO_B(H))=1$. By duality $h^0(\calO_B(K_B-H))=1$. Let
$\sigma$ be a non-zero section of $\calO_B(K_B-H)$. As usual we can write
$B=Y+Z$ where $Z$ is the maximal subcurve where $\sigma$
vanishes. Note that $Z\neq\emptyset$, since $K_B-H$ has negative
degree. Then $Y.(K_Y-H)\geq0$. By the very ampleness of $H$ this implies
$p(Y)\geq3$ and hence $p(Y)=3$. Then we must have $H.Y=4$ and by the previous
analysis $Y$ is a plane quartic with $H_Y=K_Y$.
\end{Proof}

At this point it is useful to introduce the following concept.

\begin{definition}
We say that an element $D\in |D_i|$ fulfills condition (C) if for every
decomposition $D=A+B$:

\noindent $\on{(i)}$ $p(A),p(B)\leq2$

\noindent $\on{(ii)}$ $H.A\geq 2p(A)+1$, $H.B\geq 2p(B)+1$.
\end{definition}

\begin{remark}\label{IIIrem5}
It follows immediately from (\ref{IIIgl1}) that an element $D\in |D_i|$ which
fulfills condition (C) is $3$-connected.
\end{remark}

For future use we also note

\begin{lemma}\label{IIIlemma6}
Let $D$ be a curve of genus $4$, and let $H$ be divisor on $D$ of degree $6$
with $h^0(\calO_D(H))\geq4$. Assume that for every proper subcurve $Y$ of $D$
we have $H.Y\geq 2p(Y)-1$. Then $H$ is the canonical divisor on $D$.
\end{lemma}
\begin{Proof}
By Riemann-Roch and duality $h^0(\calO_D(K_D-H))\geq1$. Let $\sigma$ be a
non-zero section of $\calO_D(K_D-H)$. As usual this defines a decomposition
$D=Y+Z$ where $Z$ is the maximal subcurve where $\sigma$ vanishes. If
$Z=\emptyset$ the claim is obvious. Otherwise $(K_D-H).Y\geq Z.Y$ and by
adjunction this gives $H.Y\leq 2p(Y)-2$, a contradiction.
\end{Proof}

Our next aim is to analyze the condition $h^0(\calO_S(H))=5$. For this
we introduce the divisor
$$
\Delta_i\equiv H-(L-x_i).
$$

\begin{lemma}\label{IIIlemma7}
The following conditions are equivalent:

\noindent $\on{(i)}$ $h^0(\calO_S(H))=5$ (resp.~$h^1(\calO_S(H))=1$).

\noindent $\on{(ii)}$ $h^0(\calO_D(H))=4$ (resp.~$h^1(\calO_D(H))=1$) for some
(every) element $D\in |D_i|$.

\noindent $\on{(iii)}$ $h^0(\calO_D(K_D-H))=1$ for some (every) element $D\in
|D_i|$.

\noindent Moreover assume that $D\in |D_i|$ fulfills condition $\on{(C)}$. Then
the following conditions are equivalent to $\on{(i)}$-$\on{(iii)}$:

\noindent $\on{(iv)}$ $\calO_D(H)=K_D$.

\noindent $\on{(v)}$ $\Delta_i|_D\equiv (2L-\sum x_i)|_D$.
\end{lemma}
\begin{Proof}
Since $h^0(\calO_S(D_i))\geq1$ we have an exact sequence
$$
0\longrightarrow\calO_S(x_i)\longrightarrow\calO_S(H)
\longrightarrow\calO_D(H) \longrightarrow 0.
$$
Since $h^1(\calO_S(x_i))=0$ the equivalence of (i) and (ii) follows. The
equivalence of (ii) and (iii) is a consequence of Serre duality. It follows
from Lemma~(\ref{IIIlemma6}) that (iii) implies (iv) if (C) holds. Conversely
if $\calO_D(H)=K_D$ then $h^0(\calO_D(K_D-H))=h^0(\calO_D)=1$, since $D$ is
$3$-connected. To show the equivalence of (iv) and (v) note that by adjunction
$$
K_D\equiv (K_S+D)|_D\equiv (3L-2x_i-\sum_{j\neq i}x_j)|_D.
$$
Hence $K_D\equiv H|_D\equiv (\Delta_i+(L-x_i))|_D$ if and only if
$\Delta_i|_D\equiv (K_D-(L-x_i))|_D\equiv (2L-\sum x_i)|_D$.
\end{Proof}

We want to discuss necessary open conditions which must be fulfilled if $|H|$
is ample.

\begin{definition}
We say that $|H|$ fulfills condition (P) if for every divisor $Y$ on $S$ with
$Y.L\leq6$, $p(Y)\leq2$, $H.Y\leq2p(Y)$ the linear system $|Y|$ is empty.
\end{definition}

\begin{remark}\label{IIIrem8}
(i) By Proposition (\ref{IIIprop3}) this condition is necessary for $|H|$
to be very ample.

\noindent (ii) Note that in order to check (P) one only need check {\em
finitely many} open conditions.

\noindent (iii) For $Y.L=0$ condition (P) implies that the only points which
can have infinitely near points are the $x_i$. The only possibility is that at
most one of the points $y_k$ is infinitely near to some point $x_i$.

\noindent (iv) If $Y.L=1$ then (P) implies
$$
|L-\sum_{i\in\triangle}x_i-\sum_{k\in\triangle'}y_k|=\emptyset \text{ for
}2|\triangle|+|\triangle'|\geq 6.
$$
In particular no three of the points $x_i$ can lie on a line.

\noindent (v) If $Y.L=6$ then (P) gives
$$
|D_i-x_j|=\emptyset\ (j\neq i),\quad |D_i-y_k-y_l|=\emptyset\ (k\neq l).
$$
\end{remark}

There are, however, two more open conditions which are not as obvious to see.

\begin{proposition}\label{IIIprop9}
If $|H|$ embeds $S$ into $\PP^4$ then the following open conditions hold:
\begin{enumerate}
\item[$\on{(Q)}$]
$|D_i-2x_i|=\emptyset, \quad
|D_i-x_i-y_k|=\emptyset, \quad
|D_i-2y_k|=\emptyset$

\item[$\on{(R)}$]
For any effective curve $C$ with $C\equiv L-x_i-x_j-y_k, C\equiv L-x_i-x_j$ or
$ C\equiv y_k$ one has $\dim |D_i-C|\le 0$. Moreover $\dim |H-(L-x_i-x_j)|\le
1$.
\end{enumerate}
\end{proposition}
\begin{Proof}
We start with (R). We already know that $\dim |D_i|=1$. Hence we have to
see that such a curve $C$ is not contained in the plane spanned by the
conic
$x_i$. But this would contradict very ampleness since $C.x_i=1$ or $0$. If
$|H|$ is very ample then it embeds $\Lambda_{ij}=L-x_i-x_j$ as a plane conic
(irreducible or reducible but reduced). The claim then follows from the exact
sequence
$$
0\longrightarrow\calO_S(H-(L-x_i-x_j))\longrightarrow \calO_S(H)
\longrightarrow \calO_{\Lambda_{ij}}(H) \longrightarrow 0.
$$

Next we consider the linear system $|D_i-2x_i|$. Assume there is a curve $B\in
|D_i-2x_i|$. Then $p(B)=-3$. Since $H.B=2$ we have the following
possibilities: $B$ is a reduced conic (either smooth or reducible). Then
$p(B)=0$, a contradiction. If $B$ is the union of $2$ skew lines, then
$p(B)=-1$ which is also not possible. Hence $B$ must be a non-reduced line.
But this is not possible, since the class of $B$ on $S$ is not divisible by
$2$.

The crucial step is to prove the
\begin{claim}
Set $D=D_i$. If $|D|$ contains $y_k+B$, then $B$ is of the form
$B=B'+(L-x_i-x_j-y_k)$ with $H_{B'}=K_{B'}$.
\end{claim}

It follows from Lemma~(\ref{IIIlemma7}) that there exists a non-zero section
$0\neq\sigma\in H^0(\calO_D(K_D-H))$. As usual this defines a decomposition
$D=Y+Z$. Since $(K_D-H).y_k=-1$ the curve $Z$ must contain the irreducible
curve $y_k$. Moreover since $y_k.B=2$ and $(K_D-H).B=1$ it follows that $Z$
contains some further curve $Z'$ contained in $B$, i.e.~$B=B'+Z'$. Now as in
proof of Lemma~(\ref{IIIlemma6}) $H.B'\leq2p(B')-2$ and very ampleness of
$|H|$ together with (\ref{IIIlemma1}) implies $p(B')=3$. As in the proof of
Proposition~(\ref{IIIprop4}) one concludes that $H.B'=4$, $H_{B'}=K_{B'}$. In
particular $Z'$ is a line. Since $p(D_i-2y_k)=1$ it follows that $Z'\neq y_k$.
First assume that $Z'.y_k=0$. Then $p(Z'+y_k)=-1$ and $B'.y_k=2$. It follows
from (\ref{IIIgl1}) that $B'.Z'=1$. But now the decomposition $Z'+(B'+y_k)$
contradicts $2$-connectedness. Hence $Z'$ and $y_k$ are two lines meeting in a
point. This gives $p(y_k+Z')=0$, $B'.(y_k+Z')=2$. We can write
$$
Z'=aL-\beta_ix_i-\sum_{j\neq i}\beta_jx_j-y_k-\sum_{l\neq k}\alpha_ly_l.
$$
If $a=0$ then $Z'=x_i-y_k$ or $Z'=x_j-y_k$, $j\neq i$. The first is impossible
since $p(D_i-x_i)=1$ the second contradicts $|D_i-x_j|=\emptyset$. Hence
$1\leq a\leq6$. Since $Z'$ is mapped to a line in $\PP^4$ we find
$Z'.y_l\leq1$, $Z'.x_j\leq2$, i.e.
\begin{gather}
0\leq\alpha_l\leq1,\quad 0\leq\beta_i,\beta_j\leq2.\label{IIIgl5}
\end{gather}
It follows from (\ref{IIIgl5}) and from $p(Z')=0$ that $a\leq4$; moreover
$p(Z')=0$, $p(B')=3$ and $p(B)=3$ imply $Z'.B'=1$. Using $0\leq\alpha_l\leq 1$
this gives
\begin{gather}
a(6-a)-\beta_i(3-\beta_i)-\sum_{j\neq i}\beta_j(2-\beta_j)=2.
\label{IIIgl6}
\end{gather}
In view of (\ref{IIIgl5}) this shows $a(6-a)\leq7$ and since $a\leq4$ it
follows that $a=1$. Then $\beta_i,\beta_j\leq1$. If $\beta_i=0$ then by
(\ref{IIIgl6}) $\beta_j=1$ for $j\neq i$, but no three of the points $x_i$ can
be collinear by (\ref{IIIgl6}). Hence $\beta_i=1$ and exactly one $\beta_j$ is
$1$. Together with $H.Z'=1$ this gives $Z'=L-x_i-x_j-y_k$ as claimed.

We are now in a position to prove that $|D_i-x_i-y_k|=\emptyset$ and
$|D_i-2y_k|=\emptyset$. For this we have to show that $B'$ cannot contain
$x_i$ or $y_k$. In the first case $B'=x_i+B''$. Then $H.x_i=2$ and
$K_{B'}.x_i=1$ contradicting $H_{B'}=K_{B'}$. Similarly in the second case
$B'=y_k+B''$ with $H.y_k=1$ and $K_{B'}.y_k=0$ giving the same contradiction.
\end{Proof}

Observe for future use that in the following proposition the assumption
that $|H|$ is very ample is not made.

\begin{proposition}\label{IIIprop10}
Assume that the open conditions $\on{(P)}$ and $\on{(Q)}$ hold. Then an
effective decomposition $D=A+B$ either fulfills condition $\on{(C)}$ and
hence is not
$3$-disconnecting or (after possibly interchanging $A$ and $B$) $A=y_k$,
$L-x_i-x_j$ or $L-x_i-x_j-y_k$.
\end{proposition}
\begin{Proof}
Let $D=A+B$. Clearly we can assume $A.L\leq3$. We shall first treat the case
$A.L=0$, i.e.~$A$ is exceptional with respect to the blowing down map
$S\to\PP^2$. Then $p(A)\leq0$ and $A.H>0$ by (P). By conditions (Q) and
(P) (cf. Remark (III.8)(v)) if $A.H=1$, then either $A=x_j-y_k$ or
$A=x_i-y_k$ or $A=y_k$. In the first two cases $A.B\ge 3$ and $p(B)\le
2$, the third is one of the exceptions stated. If $A.H\ge 2$ then
$p(B)\le 2$ and the claim follows from (P).

Hence we can now write
\begin{eqnarray*}
A&\equiv&aL-\sum\alpha_jx_j-\sum a_ky_k\\
B&\equiv&bL-\sum\beta_jx_j-\sum b_ky_k
\end{eqnarray*}
with $a,b>0$. Using the open conditions from Remark~(\ref{IIIrem8})(v) (which
are a consequence of (P))  and (Q) it follows that
\begin{gather*}
\begin{aligned}
a_k,b_k &\geq -1,\\
\alpha_j,\beta_j &\geq 0,\\
\alpha_i,\beta_i &\geq -1,
\end{aligned}
\quad
\begin{aligned}
a_k+b_k &= 1\\
\alpha_j+\beta_j &= 2\\
\alpha_i+\beta_i &= 3
\end{aligned}
\quad
(j\neq i)
\end{gather*}
and moreover that at most one of the integers $a_k,b_k,\alpha_i,\beta_i$ can
be negative. If $\beta_i=-1$ then $\alpha_i=4$. In this case $A$ cannot be
effective since we have assumed $a\leq3$. If $\alpha_i=-1$ then $\beta_i=4$
and hence $b\geq4$. We have to consider the cases $a=1$ or $2$. In either case
$p(A)\leq0$ and $H.A\geq2p(A)+1$ follows from (P). On the other hand
\begin{eqnarray*}
H.B-(2p(B)+1)&=&(9b-b^2+1)+\sum_{j\neq i}\beta_j(\beta_j-3)+\sum_k b_k(b_k-2)\\
&\geq&(9b-b^2+1)-6-12\geq3
\end{eqnarray*}
since $b=4,5$. Hence we can now assume $\alpha_i,\beta_i\geq0$.

\noindent $\mathbf{a=1}$. We first treat the case $a_k\geq0$ for all $k$. Then
$$
A\equiv L-\sum_{j\in\triangle}x_j-\sum_{k\in\triangle'} y_k.
$$
Clearly $p(A)\leq0$. Let $\delta_{i\triangle}=0$ (resp.~$1$) if
$i\not\in\triangle$ (resp.~$i\in\triangle$). Then
$$
p(B)=|\triangle|+\delta_{i\triangle}.
$$
We only have to treat the cases where $p(B)\geq3$. Then either
$\delta_{i\triangle}=0$, $|\triangle|\geq3$ or $\delta_{i\triangle}=1$,
$|\triangle|\geq2$. In the first case
$$
H.A=6-2|\triangle|-|\triangle'|\leq0
$$
contradicting (P) for $A$. In the second case the only possibilty is
$|\triangle|=2$, $|\triangle'|\leq1$. But then $A=L-x_i-x_j$ or
$L-x_i-x_j-y_k$. Now assume that one $a_k$ is negative. We can assume
$a_{16}=-1$. Then
$$
A\equiv L-\sum_{j\in\triangle}x_j-\sum_{k\in\triangle'}y_k+y_{16}.
$$
In this case $p(A)=-1$ and
$$
p(B)=|\triangle|+\delta_{i\triangle}-1.
$$
Using the same arguments as before we find that $p(B)\leq2$ in all cases.

\noindent $\mathbf{a=2}$. Again we first assume that all $a_k\geq0$. Then
$$
A\equiv 2L-\sum_{j\in\triangle}x_j-\sum_{k\in\triangle'}2x_k-
\sum_{l\in\triangle''}y_l-\sum_{m\in\triangle'''}2y_m.
$$
Clearly $p(A)\leq0$. If $i\not\in\triangle\cup\triangle'$ then $p(B)\leq0$. If
$i\in\triangle$ then $p(B)\leq2$. Now assume that $i\in\triangle'$. In this
case $p(B)\leq2$ with one possible exception: $|\triangle|=3$ and
$|\triangle'''|=0$. But then
$$
A\equiv 2L-2x_i-x_j-x_k-x_l-\sum_{l\in\triangle''}y_l.
$$
In this case $A$ splits into two lines meeting $x_i$. But then one of these
lines must contain $3$ of the points $x_j$ contradicting condition (P).
Finally let $a_{16}=-1$. The above arguments show that in this case
$p(B)\leq2$.

\noindent $\mathbf{a=3}$. Since in this case $p(A),p(B)\leq1$ condition
(C) follows.
\end{Proof}

Propositions~(\ref{IIIprop4}) and (\ref{IIIprop10}) have provided us with a
fairly good understanding of the behaviour of $H$ on the pencil $|D_i|$.

\begin{cor}\label{IIIcor11}
Assume $|H|$ embeds $S$ into $\PP^4$. For every element $D\in |D_i|$ either:

\noindent $\on{(i)}$ $D$ is $3$-connected and $H_D=K_D$ or

\noindent $\on{(ii)}$ $D=B+(L-x_i-x_j)$ with $H_B=K_B$.
\end{cor}

\begin{remark}\label{IIIrem12}
The conic $L-x_i-x_j$ can be irreducible or reducible in which case it splits
as $(L-x_i-x_j-y_k)+y_k$.
\end{remark}

At this point we can also conclude our discussion about the linear system
$|\Delta_i|=|H-(L-x_i)|$ (cf.~(\ref{IIIlemma7})).

\begin{proposition}\label{IIIprop13}
If $|H|$ embeds $S$ into $\PP^4$, then $\dim|\Delta_i|=0$.
\end{proposition}
\begin{Proof}
We first claim that the general element $D\in |D_i|$ is $3$-connected. Indeed
if $D$ is not $3$-connected, then $D=B+(L-x_i-x_j)$. The conic $L-x_i-x_j$
spans a plane $E'$. If $E$ is the plane spanned by $x_i$ then $E\neq E'$
since $(L-x_i-x_j).x_i=1$. Hence $D$ is cut out by the hyperplane spanned by
$E$ and $E'$. Varying the index $j$ there are at most $3$ such hyperplanes.

Clearly $L-x_i$ is effective. Consider the exact sequence
$$
0\longrightarrow\calO_S(\Delta_i)\longrightarrow\calO_S(H) \longrightarrow
\calO_S(H)|_{L-x_i} \longrightarrow 0.
$$
Since $H.(L-x_i)=4$ and $p(L-x_i)=0$ it follows that $|H|$ cannot map $L-x_i$
to a plane curve. This shows $h^0(\calO_S(\Delta_i))\leq1$.

On the other hand choose an element $D\in|D_i|$ which is $3$-connected. We
have an exact sequence
$$
0\longrightarrow\calO_S(2x_i-L)\longrightarrow\calO_S(\Delta_i)
\longrightarrow \calO_D(\Delta_i) \longrightarrow 0.
$$
Now $h^0(\calO_S(2x_i-L))=h^2(\calO_S(2x_i-L))=0$ and hence
$h^1(\calO_S(2x_i-L))=1$ by Riemann-Roch. Since $|H|$ is ample no $3$ of the
points $x_i$ lie on a line. Hence $|2L-\sum x_i|$ is a base point free pencil.
Since $|(2L-\sum x_i)-D|=\emptyset$ this shows that $|2L-\sum x_i|$ cuts out a
base-point free pencil on $D$. Since $D$ is $3$-connected $(2L-\sum
x_i)|_D\equiv\Delta_i|_D$ by Lemma~(\ref{IIIlemma7}) and hence
$h^0(\calO_D(\Delta_i))\geq2$. By the above sequence this implies
$h^0(\calO_S(\Delta_i))\geq 1$.
\end{Proof}

We are now ready to characterize very ample linear systems which embed $S$
into $\PP^4$.

\begin{theorem}\label{IIItheo14}
The linear system $|H|$ embeds $S$ into $\PP^4$ if and only if

\noindent $\on{(i)}$ The open conditions $\on{(P)}$, $\on{(Q)}$ and $\on{(R)}$
hold.

\noindent $\on{(ii)}$ The following closed conditions hold:
\begin{enumerate}
\item[$\on{(}D_i\on{)}$] $\dim |D_i|=1$
\item[$\on{(}\Delta_i\on{)}$] For a $3$-connected element $D\in|D_i|$ (whose
existence follows from the above conditions) $\Delta_i.D\equiv (2L-\sum
x_i).D$.
\end{enumerate}
\end{theorem}

\begin{remark}\label{IIIrem15}
As the proof will show it is enough to check the closed conditions ($D_i$),
($\Delta_i$) for one $i$.
\end{remark}

\begin{Proof}
We have already seen that these conditions are necessary. Next we shall show
that a $3$-connected element $D\in|D_i|$ exists if the open conditions and
($D_i$) are fulfilled. Assume that no element $D\in|D_i|$ is $3$-connected.
Then by Proposition~(\ref{IIIprop10}) every element $D$ is of the form $D=B+C$
with $C=L-x_i-x_j$, $L-x_i-x_j-y_k$ or $y_k$. But by condition (R) there are
only finitely many such elements in $|D_i|$.

We shall now proceed in several steps.

\noindent {\bf Step 1}: $h^0(\calO_S(H))=5$.

We have seen in the proof of Lemma~(\ref{IIIlemma7}) that for a $3$-connected
element $D$ the equality $\Delta_i.D\equiv (2L-\sum x_i).D$ implies
$K_D=H_D$ and hence $h^0(\calO_D(K_D-H))=1$, resp.~$h^1(\calO_D(H))=1$. Now the
claim follows from the equivalence of (i) and (ii) in Lemma~(\ref{IIIlemma7}).

In order to prove very ampleness of $|H|$ we want to apply the Alexander-Bauer
Lemma to the decomposition
$$
H\equiv D_i+x_i.
$$
We first have to show that $|H|$ cuts out complete linear systems on $x_i$ and
$D\in|D_i|$. Recall that $x_i$ is either a $\PP^1$ or consists of two
$\PP^1$'s meeting transversally (cf.~Remark~(\ref{IIIrem8})(iii)). Moreover
$H.x_i=2$ and if $x_i$ is reducible then $H$ has degree $1$ on every
component. Hence
$h^0(\calO_{x_i}(H))=3$. The claim for $x_i$ then follows from the exact
sequence
$$
0\longrightarrow\calO_S(D_i)\longrightarrow\calO_S(H) \longrightarrow
\calO_{x_i}(H)\longrightarrow 0.
$$
and condition ($D_i$), i.e.~$h^0(\calO_S(D_i))=2$. The corresponding claim for
$D$ follows from the sequence
$$
0\longrightarrow\calO_S(x_i)\longrightarrow\calO_S(H) \longrightarrow
\calO_S(H)|_D \longrightarrow 0.
$$
Our above discussion also shows that $|H|$ embeds $x_i$ as a conic (which can
be irreducible or consist of two different lines).

\noindent {\bf Step 2}: If $D\in|D_i|$ is $3$-connected then $H_D=K_D$ and
$|H|$ is very ample on $D$.

We have already seen the first claim. We have to see that $K_D$ is very ample.
For this we consider the pencils $|\Sigma_1|=|L-x_i|$,
resp.~$|\Sigma_2|=|2L-\sum x_j|$. Clearly $|\Sigma_1|$ is base point free and
the same is true for $|\Sigma_2|$ as no three of the points $x_i$ lie on a
line (by (P)). Hence
$$
|\Sigma_1+\Sigma_2|=|3L-2x_i-\sum_{j\neq i}x_j|=|D_i+K_S|
$$
is base point free. By adjunction $(D_i+K_S)|_D\equiv K_D$ and the exact
sequence
$$
0\longrightarrow\calO_S(K_S)\longrightarrow\calO_S(K_S+D_i) \longrightarrow
\calO_D(K_D) \longrightarrow 0
$$
shows that restriction defines an isomorphism $|\Sigma_1+\Sigma_2|\cong
|K_D|$. Let $X$ be the blow-up of $\PP^2$ in the points $x_j$ and $\pi:S\to X$
the map blowing down the exceptional curves $y_k$. The linear system
$|\Sigma_1+\Sigma_2|$ defines a morphism
$$
f=\phi_{|\Sigma_1+\Sigma_2|}:X\longrightarrow\PP^3.
$$
It is easy to understand the map $f$: Clearly $f$ contracts the three
$(-1)$-curves $\Lambda_{ij}=L-x_i-x_j$, $j\neq i$. Let $\pi':X\to X'$ be the
map which blows down the curves $\Lambda_{ij}$ (this makes also sense if
$\Lambda_{ij}=(L-x_i-x_j-y_k)+y_k$ where we first contract $y_k$ and then
$L-x_i-x_j-y_k$). Then $X'$ is a smooth surface and we have a commutative
diagram
$$
\unitlength1cm
\begin{picture}(4,2)
\put(0,1.5){$X$} \put(0.5,1.6){\vector(1,0){2.4}} \put(3.1,1.5){$f(X)$}
\put(0.5,1.3){\vector(3,-2){1}} \put(1.5,0.1){$X'$}
\put(2,0.62){\vector(3,2){1}}
\put(0.5,0.7){{$\scriptstyle \pi'$}}
\put(2.5,0.7){{$\scriptstyle f'$}}
\put(1.5,1.7){{$\scriptstyle f$}}
\end{picture}
$$
where $f'$ maps $X'$ isomorphically onto a smooth quadric. This shows that
$\phi_{|K_D|}:D\to\PP^3$ is the composition of the blowing down maps $\pi:S\to
X$ and $\pi':X\to X'=\PP^1\times\PP^1$ followed by an embedding of $X'$. Now
$D.y_k=1$, hence $\pi|_D$ can only fail to be an isomorphism if $D$ contains
$y_k$. But this is impossible if $D$ is $3$-connected. Similarly
$D.\Lambda_{ij}=1$ and $D$ cannot contain a component of $\Lambda_{ij}$. Hence
we are done in this case.

It remains to treat the case when $D$ is not $3$-connected.

\noindent {\bf Step 3}: If $D$ is not $3$-connected, then $D=B+(L-x_i-x_j)$,
$H_B=K_B$ and $|H|$ restricts onto $|K_B|$.

We have already seen that $h^0(\calO_S(H))=5$ and hence
$h^0(\calO_D(K_D-H))=1$. As usual a non-zero section $\sigma$ defines a
decomposition $D=Y+Z$. Our first claim is that $Z$ is different from $0$.
In fact if
$Z=0$ then
$K_D-H$ would be trivial on
$D$. On the other hand $D$ is not $3$-connected, thus it splits as $D=A+B$
with
$A$ as in Proposition~(\ref{IIIprop10}), in particular $p(A)=0$, $A.B=2$.
Then $K_D.A=0$ contradicting $H.A>0$ which follows from (P). Thus $Z$ is
different from $0$ and since the section $\sigma$ defines a good section
$\sigma'$ of
$H^0({\cal O}_Y(K_Y-H))$ it follows that $2p(Y)-2\ge H.Y$, and hence
$p(Y)\ge 3, Y.Z\le 2$. Then Proposition (\ref{IIIprop10}) applies and
$Z=y_k$ or
$L-x_i-x_j-y_k$ or $L-x_i-x_j$.
If $Z=y_k$ or
$L-x_i-x_j$ then
$(K_Y-H).Y=-1$, a contradiction. Hence $Z=L-x_i-x_j$ and $H_Y=K_Y$. We next
claim that $B$ is $2$-connected. Assume we have a decomposition $B=B_1+B_2$
with $B_1.B_2\leq1$. Then $(B_1+B_2).(L-x_i-x_j)=2$, hence we can assume that
$B_1.(L-x_i-x_j)\leq1$. But then $B_1.(B_2+L-x_i-x_j)\leq2$ contradicting
Proposition~(\ref{IIIprop10}). This shows that $h^1(\calO_B(K_B))=1$ and
$h^0(\calO_B(K_B))=3$. The claim then follows from the exact sequence
$$
0\longrightarrow\calO_S(L-x_j)\longrightarrow\calO_S(H) \longrightarrow
\calO_B(H) \longrightarrow 0.
$$

\noindent {\bf Step 4}: $|H|$ embeds $D$.

Our first claim is that $|H|$ embeds $B$ as a plane quartic. Since $B-y_k$
is not effective by condition (P) and $B.y_k=1$ it follows that the curve
$B$ is mapped isomorphically onto its image under the blowing down map
$\pi:S\to X$. On $X$
$$
B\equiv 5L-2x_i-x_j-2x_k-2x_l,\quad K_B\equiv (2L-x_i-x_k-x_l)|_B.
$$
Thus $|K_B|$ is induced by a standard Cremona transformation centered at
$x_j$, $x_k$ and $x_l$. Again by (P) it follows that $B-\Lambda_{ik}$ for
$k\neq i$ and $B-\Lambda_{kl}$ for $k,l\neq i$ are not effective. Since
$B.\Lambda_{ik}=B.\Lambda_{kl}=1$ it follows that $B$ is mapped isomorphically
onto a plane quartic.

It follows from condition (R) that $|H|$ embeds $\Lambda_{ij}$ as a plane conic
$Q$. The planes containing $B$ and $Q$ intersect in a line and span a $\PP^3$.
The line of intersection cannot be a component of $Q$ since, by taking
residual intersection with hyperplanes containing $B$, this would contradict
$h^0(\calO_S(x_i+y_k))=1$, resp.~$h^0(\calO_S(L-x_j-y_k))=1$. Hence the
schematic intersection of the embedded quartic $B$ and the conic $Q$ has
length at most $2$. Let $D'$ be the schematic image of $D$. Then $\calO_{D'}$
is contained in the direct image of $\calO_D$. But the former has colength
$\leq2$ in $\calO_Q\oplus\calO_B$, the latter has colength $2$, thus $D=D'$.
\end{Proof}

\begin{remark}\label{IIIrem16}
We have already remarked that conditions (P) and (Q) lead to finitely many
open conditions. Going through the proof of
Proposition~(\ref{IIIprop10}) one sees that it is sufficient to check that no
decomposition $A+B=D\in|D_i|$ exists where $A$ (or $B$) contradicts one of the
following conditions below: Here $\triangle$ and $\triangle'$ are always
disjoint subsets of
$\{1,\ldots,4\}$ whereas $\triangle''$ is a subset of $\{5,\ldots,16\}$. We
set $\delta_{i\triangle}=1$ (resp.~$0$) if $i\in\triangle$
(resp.~$i\not\in\triangle$). Similarly we define $\delta_{i\triangle'}$.
Moreover $\delta_{m}=1$ for at most one $m\in\{5,\ldots,16\}$ and
$\delta_{m}=0$ otherwise. If $\delta_{m}=1$ then ~$m\not\in\triangle''$.

\begin{enumerate}
\item[(0)]
$|x_j-x_k|=\emptyset$ ($j\neq k$), $|y_k-y_l|=\emptyset$
($k\neq l$), $|y_k-x_j|=\emptyset$, $|x_j-y_k-y_l|=\emptyset$.

\item[(1)]
$|L-\sum\limits_{j\in\triangle}x_j-
\sum\limits_{k\in\triangle''}y_k|=\emptyset$ for
$2|\triangle|+|\triangle'|\geq6$

\item[(2)]
$|2L-\sum\limits_{j\in\triangle}x_j-
\sum\limits_{k\in\triangle''}y_k|=\emptyset$ for
$2|\triangle|+|\triangle'|\geq12$.

\item[(3)]
$|3L-2x_j-\sum\limits_{k\in\triangle}x_k-
\sum\limits_{l\in\triangle''}y_l|=\emptyset$ for
$2|\triangle|+|\triangle''|\geq14$\\
 $|3L-\sum\limits_{j\in\triangle}x_j-2y_k-
\sum\limits_{l\in\triangle''}y_l|=\emptyset$ for
$2|\triangle|+|\triangle''|\geq16$\\
 $|3L-\sum\limits_{j\in\triangle}x_j-
\sum\limits_{k\in\triangle''}y_k|=\emptyset$ for
$2|\triangle|+|\triangle''|\geq16$

\item[(4)]
$|4L-(3-\delta_{i\triangle}-2\delta_{i\triangle'})x_i-
\sum\limits_{j\neq i \atop j\in\triangle}x_j- 2\sum\limits_{k\neq i \atop
k\not\in(\triangle\cup\triangle')}x_k-
\sum\limits_{l\not\in\triangle''}y_l-\delta_my_m|=\emptyset$ for
 $|\triangle|+|\triangle'|+\delta_{i\triangle}+
2\delta_{i\triangle'}-\delta_m\leq5$,
$2|\triangle'|+|\triangle''|-2\delta_{i\triangle}-4\delta_{i\triangle'}+
\delta_m\leq0$, $2|\triangle|+4|\triangle'|+|\triangle''|\leq11$

\item[(5)]
$|5L-(3-\delta_{i\triangle})x_i-\sum\limits_{j\neq i \atop
j\in\triangle}x_j- 2\sum\limits_{k\neq i \atop k\not\in\triangle}x_k-
\sum\limits_{l\not\in\triangle''}y_l-\delta_my_m|=\emptyset$ for
 $|\triangle|+\delta_{i\triangle}-\delta_m\leq2$,
$|\triangle''|-2\delta_{i\triangle}+\delta_m\leq0$,
$2|\triangle|+|\triangle''|\leq5$.

\item[(6)]
$|D_i-x_j|=\emptyset$ ($i\neq j$), $|D_i-2x_i|=\emptyset$,
$|D_i-x_i-y_k|=\emptyset$, $|D_i-2y_k|=\emptyset$, $|D_i-y_k-y_l|=\emptyset$
($k\neq l$).
\end{enumerate}
\end{remark}

Now we want to show how Theorem (\ref{IIItheo14}) can be used to prove the
existence of the special surfaces of degree 8 by explicitly constructing a
very ample linear system $|H|$. Let $x_1,\ldots , x_4$ be points in general
position in $\PP^2$, and blow them up. The linear system
$|5L-x_1-2\sum\limits_{j\ge 2} x_j|$ is 10-dimensional, its elements have
arithmetic genus 3. Let $\Delta_1$ be a general (and hence smooth) element of
the 10-dimensional linear system $|5L-x_1-2\sum\limits_{j\ge 2} x_j|$ on
${\hat\PP}^2=\PP^2(x_1, \ldots , x_4)$. Note that the image of $\Delta_1$ in
$\PP^2$ is the image of the canonical model of $\Delta_1$ under a standard
Cremona transformation. The linear system $|2L-\sum\limits_j x_j|$ cuts out a
$g^1_3$ on $\Delta_1$, since $H^1({\hat \PP}^2, {\cal O}_{{\hat
\PP}^2}(-3L+\sum x_j))=0$. The linear system
$$
|L_0|:=|(6L-3x_1-2\sum\limits_{j\ge 2}
x_j)|_{\Delta_1}-g^1_3|=|(4L-2x_1-\sum\limits_{j\ge 2} x_j)|_{\Delta_1}|
$$
on $\Delta_1$ has degree 12 and dimension 9. The linear system
$|4L-2x_1-\sum\limits_{j\ge 2} x_j|$ on ${\hat\PP}^2$ cuts out a subsystem of
codimension 1 in $|L_0|$. We consider the variety
$$
{\cal M}:=\{(\Delta_1, \sum y_k);\ \Delta_1 \mbox{ smooth }, \sum y_k \in
|L_0|\}.
$$
${\cal M}$ is rational of dimension 19.

\begin{theorem}\label{IIItheo17}
There is a non-empty open set ${\cal U}$ of the rational variety ${\cal M}$ for
which the linear system $|H|$ embeds $S$ into $\PP^4$.
\end{theorem}
\begin{Proof}
We have to show that for a general choice of $\Delta_1$ and $\sum y_k \in
|L_0|$ the linear system $|H|$ embeds $S$ into $\PP^4$. We shall first treat
the closed conditions. Since $\Delta_1$ is smooth we can identify it with its
strict transform on $S$. Consider the exact sequence
$$
0\rightarrow {\cal O}_S(L-2x_1)\rightarrow {\cal
O}_S(D_1)\rightarrow  {\cal O}_{\Delta_1}(D_1)\rightarrow 0.
$$
Since $\sum y_k\in |L_0|$ we have
$$
(35) \quad 6L-3x_1-2\sum\limits_{j\ge 2} x_j-\sum y_k\equiv g^1_3 \mbox{ on
} \Delta_1
$$
and hence $h^0({\cal O}_S(D_1))=h^0({\cal O}_{\Delta_1} (D_1))=2$. This is
condition $\on {(D_1)}$. Condition $ (\Delta_1)$ holds by construction.\\
In order to treat the open conditions we will first consider special points in
${\cal M}$ which give us all open conditions but two. These we will then treat
afterwards. The linear system $|4L-2x_1-\sum\limits_{j\ge 2} x_j|$ is free on
${\hat \PP}^2$. Hence a general element $\Gamma$ is smooth and intersects
$\Delta_1$ transversally in 12 points $y_k$ which neither lie on an exceptional
line, nor on a line of the form $\Lambda_{kl}=L-x_k-x_l$. Moreover a general
element
$\Gamma$ is irreducible. This follows since $\Gamma^2=9$ and $|\Gamma|$ is
not composed of a pencil, since the class of $\Gamma$ is not divisible by 3
on ${\hat \PP}^2$. Let $\Gamma'$ be the smooth transform of $\Gamma$ on $S$.
Since $\Gamma$ is smooth, $\Gamma'$ is isomorphic to $\Gamma$.
\begin{claim}
$|D_1|=\Gamma'+|2L-\sum\limits_j x_j|.$\\
This follows immediately since $D_1\equiv \Gamma'+(2L-\sum\limits_j x_j)$ and
$\dim |D_1|=1=\dim(\Gamma'+|2L-\sum x_j|).
$
\end{claim}
The only curves contained in an element of $|D_1|$ are $\Gamma'$, conics
$C\equiv 2L-\sum x_j$ and lines $\Lambda_{kl}=L-x_k-x_l$. The latter only
happens for finitely many elements of $|D_1|$. This shows immediately that
conditions (Q) and (R) are fulfilled with the possible exception that
$\dim|H-\Lambda_{1j}|\ge 2$. To exclude this we consider w.l.o.g. the
case $j=2$. Note that $H-\Lambda_{12}\equiv \Delta_2+x_1\equiv\Gamma ' +
\Lambda_{34}+x_1$. Since $\Gamma '$ is smooth of genus 2 and  $\Gamma
'.(\Delta_2+x_1)=1$ it follows that $h^0({\cal O}_{\Gamma
'}(\Delta_2+x_1))\le 1$. The claim now follows from the exact sequence
$$
0\rightarrow {\cal O}_S(\Lambda_{34}+x_1)\rightarrow {\cal
O}_S(\Delta_2+x_1)\rightarrow {\cal O}_{\Gamma
'}(\Delta_2+x_1)\rightarrow 0
$$
together with the fact that $h^0({\cal
O}_S(\Lambda_{34}+x_1))=1$.
It remains to consider (P). The curve
$\Gamma'$ contradicts condition (P) since
$p(\Gamma')=2, H.\Gamma'=4$. Similarly the decomposition
$(\Gamma'+\Lambda_{ij})+\Lambda_{kl}$ contradicts (P) if $k, l\neq 1$. On
the other hand the above construction shows that for one (and hence the
general) pair $(\Delta_k, \sum y_k)$ all open conditions given by (P) are
fulfilled for a decomposition $D=A+B$ of an element in
$|D_1|$ with the possible exception of $|\Gamma'|\neq\emptyset$ or
$|D_1-\Lambda_{kl}|\neq\emptyset$ for $k,l\neq 1$. The first case is easy,
we can simply take an element $\sum y_k \in |L_0|$ which is not in the
codimension 1 linear subsystem given by
$|4L-2x_1-\sum\limits_{j\ge2} x_j|$ on ${\hat\PP}^2$. Next we assume that
there is an element $A\in |D_1-\Lambda_{kl}|$ where $k,l\neq 1$. Then
$A.\Delta_1=2$. Since
$\Delta_1$ cannot be a component of $A$ this means that $A$ intersects
$\Delta_1$ in two points $Q_0, Q_1$. If $j$ is the remaining element of the
set $\{1,\ldots,4\}$ then $L-x_1-x_j\equiv Q_0+Q_1$ on $\Delta_1$ The linear
system $|L|$ cuts out a $g^2_5$ on $\Delta_1$ and is hence complete. Hence
$Q_0 + Q_1$ is the intersection of $\Lambda_{1j}$ with $\Delta_1$. In
particular $\Lambda_{1j}$ intersects $A$ in at least 2 points, namely $Q_0$
and $Q_1$. Since $A.\Lambda_{1j}=0$ this implies that $\Lambda_{1j}$ is a
component of $A$ (we can assume that $\Lambda_{1j}$ is irreducible). Hence
$A=A'+\Lambda_{1j}$ with $A'\in |D_1-\Lambda_{kl}-\Lambda_{1j}|=|\Gamma'|$
and we are reduced to the previous case.
\end{Proof}

\begin{remarks}\label{IIIrem18}
(i) Originally Okonek \cite{O2} constructed surfaces of degree $8$ and
sectional genus $6$ with the help of reflexive sheaves.

\noindent (ii) According to \cite{DES} the rational surfaces of degree $8$
with $\pi=6$ arise as the locus where a general morphism
$\phi:\Omega^3(3)\to\calO(1)\oplus 4\calO$ drops rank by $1$. The space of
such maps has dimension $80$. Taking the obvious group actions into account we
find that the moduli space has dimension $43=19+\dim\on{Aut}\PP^4$.
Moreover this description shows that the moduli space is irreducible and
unirational.

\noindent (iii) These surfaces are in $(3,4)$-liaison with the Veronese
surface \cite{O2}. Counting parameters one finds again that they depend on $19$
parameters (modulo $\on{Aut}(\PP^4)$).

\noindent (iv) It was pointed out to us by K.~Ranestad that Ellingsrud and
Peskine (unpublished) also suggested a construction of these surfaces via
linear systems. They start with a smooth quartic $K_4=\{f_4=0\}$ and a
smooth quintic
$K_5=\{f_5=0\}$ touching in $4$ points $x_1,\ldots,x_4$. Let
$y_5,\ldots,y_{16}$ be the remaining points of intersection. Let
$$
\calI'=\calO_{\PP^2}\left(-\sum x_i\right),\quad
\calI=\calO_{\PP^2}\left(-2\sum x_i-\sum y_k\right).
$$
Then we have an exact sequence
$$
0\longrightarrow \calI'(-4)\longrightarrow \calI\longrightarrow \calO_{K_4}(-5)
\longrightarrow 0.
$$
Twisting this by $\calO(6)$ and taking global section gives
$$
0\longrightarrow \Gamma(\calI'(2)) \longrightarrow \Gamma(\calO_S(H))
\longrightarrow \Gamma(\calO_{K_4}(1)) \longrightarrow 0.
$$
Since $h^0(\calI'(2))=2$ and $h^0(\calO_{K_4}(1))=3$ this shows
$h^0(\calO_S(H))=5$. One can easily see that $|\Delta_i|\neq\emptyset$ and
$\dim|D_i|\geq1$ in this construction: counting parameters one shows that
$\Delta_i=\{lf_4+f_5=0\}$ for some suitable linear form and that there is
at least a $1$-dimensional family of curves in $|D_i|$ which are of the form
$D=\{qf_4+lf_5\}$ where $q$ is of degree $2$ and $l$ is a linear form. This
construction, too, depends on $19$ parameters.
\end{remarks}

Finally we want to discuss the moduli space of smooth special surfaces of
degree 8 in $\PP^4$ (modulo Aut $\PP^4$). Recall the set ${\cal M}$ consisting
of pairs $(\Delta_1, \sum y_k)$ where $\Delta_1 \in |H-(L-x_1)|$ is smooth
and $\sum y_k \in |L_0|$. We have proved in Theorem (\ref{IIItheo17}) that for
a general pair $(\Delta_1, \sum y_k)$ the linear system $|H|$ embeds $S$ into
$\PP^4$. Indeed in this way we obtain the general smooth surface of degree
8 in $ \PP^4$. The surface $X={\hat \PP}^2$, i.e. $\PP^2$ blown up in
$x_1,\ldots, x_4$ is the del Pezzo surface of degree 5. It is well known that
Aut$X\cong S_5$ the symmetric group in 5 letters (Aut $X$ acts transitively on
the 5 maximal sets of disjoint rational curves on $X$, see \cite [Chapter
IV]{M}).

\begin{proposition}\label{IIIprop19}
For general $S$ the only lines contained in $S$ are the $y_k$'s.
\end{proposition}
\begin{Proof}
Let $l$ be a line on S. The statement is clear if $l$ is $\pi$-exceptional as
the $x_i$ are mapped to conics and since we can assume that there are no
infinitesimally near points. If $l$ is not skew to the plane spanned by $x_i$
then $l$ is contained in a reducible member of $|D_i|$. But for general
choice there is no decomposition A+B with A (or B) a line. Hence we can assume
that $l.x_i=0$ for $i=1,\ldots, 4$ and $l.y_k\le 1$ for all $k$. Thus $l\equiv
a L-\sum\limits_{k\in \triangle} y_k$ with $a\le 2$. Since $H.l=1$ we have
either
$a=1$ and $|\triangle|=5$ or $a=2$ and $|\triangle|=11$. In the first case 5 of
the
$y_k$ are collinear. But then it follows with the monodromy argument of
\cite [p.111]{ACGH} that all the $y_k$'s are collinear which is absurd. In
the same way the case $a=2$ would imply that all the $y_k$'s are on a conic
which also contradicts very ampleness of $|H|$.
\end{Proof}
\begin{theorem}\label{IIItheo20}
The moduli space of polarized rational surfaces (S,H) where $|H|$ embeds
$S$ into $\PP^4$ as a surface of degree 8, speciality 1 and sectional
genus 6 is birationally equivalent to
${\cal M}/S_5$.
\end{theorem}
\begin{Proof}
Let ${\cal V}$ be the open set of ${\cal M}$ where $|H|$ embeds $S$ into
$\PP^4$ and where all the $\Delta_i$'s are smooth. Let $(\Delta_1, \sum
y_k)$ and  $(\Delta_1 ', \sum y_k ')$ be two elements which give rise to
surfaces $S, S'\subset \PP^4$ for which a projective transformation
${\bar g}:S\rightarrow S'$ exists. Since obviously  ${\bar g}$ carries
lines to lines, it follows from Proposition (\ref{IIIprop19}) that ${\bar
g}$ is induced by an automorphism $g:X\rightarrow X$ carrying the set
$\{y_k\}$ to $\{y_k'\}$. Conversely, the group $S_5=\mbox { Aut } (X)$
acts on ${\cal V}$ as follows. Let $S$ correspond to $(\Delta_1 \sum
y_k)$ and let $g\in \mbox { Aut } (X)$: Then, since $6L-2\sum
x_j=-2K_X$ which is invariant under the action of $S_5$, we set
$\{y_k'\}=g\{y_k\}, H'=-2K_X-\sum y_k'$. Then $H'$ embeds $S'={\tilde
X}(y_1', \ldots, y_{12}')$ and we set $\Delta_1'$ to be the unique
curve in $|H'-L+y_1|$.
\end{Proof}

\section{Further outlook}\label{sectionIV}

In this section we want to discuss how this method can possibly be applied to
other surfaces. For smooth surfaces of degree $\leq8$ it is rather
straightforward to give a decomposition $H\equiv C+D$ which allows to apply the
Alexander-Bauer lemma. This was done in \cite{B}, \cite{CF} and
section~\ref{sectionIII} of this article. In degree $9$ there is one
non-special surface, which was treated in section~\ref{sectionII} of this
article, and a special surface with sectional genus $\pi=7$ which was found by
Alexander \cite{A2}. Here $S$ is $\PP^2$ blown up in 15 points
$x_1,\ldots,x_{15}$ and $H\equiv 9L-3\sum\limits_{i=1}^6x_i-
2\sum\limits_{j=7}^9x_j-\sum\limits_{k=10}^{15}x_k$. As pointed out by
Alexander one can take the decomposition $H\equiv C+D$ where $C\equiv
3L-\sum\limits_{i=1}^9x_i$ and $D\equiv H-C$. Then $C$ is a plane cubic and
$|D|$ is a pencil of canonical curves of genus 4.

Rational surfaces of degree 10 were treated by Ranestad \cite{R1}, \cite{R2},
Popescu and Ranestad \cite{PR} and Alexander \cite{A2}. There is one surface
with $\pi=8$. In this case $S$ is $\PP^2$ blown up in 13 points and $H\equiv
14L-6x_1-4\sum\limits_{i=2}^{10}x_i-2x_{11}-x_{12}-x_{13}$. Following
Alexander \cite{A2} the curve $C\equiv
7L-3x_1-2\sum\limits_{i=2}^{10}x_i-\sum\limits_{j=11}^{13}x_j$ is a plane
quartic and the residual pencil $|D|$ has $p(D)=3$ and degree 6. For sectional
genus $\pi=9$ there are two possibilities. The first is $\PP^2$ blown in 18
points with $H\equiv 8L-2\sum\limits_{i=1}^{12}x_i-
\sum\limits_{j=13}^{18}x_j$. One can take $C\equiv 4L-\sum\limits_{i=1}^{16}
x_i$ which becomes a plane quartic. For the residual intersection $|D|$ one
finds $p(D)=3$, $H.D=6$. (For more details of this geometrically interesting
situation see \cite[Proposition~2.2]{PR}. The second  surface with $\pi=9$ is
more difficult. Again we have $\PP^2$ blown up in 18 points, but this time
$H\equiv 9L-3\sum\limits_{i=1}^4x_i- 2\sum\limits_{j=5}^{11}x_j-
\sum\limits_{k=12}^{18}x_k$. Clearly $S$ contains plane curves, e.g.~the
conics $x_j$. But then for the residual pencil $|D|$ one has $p(D)=7$, $H.D=9$
and this case seems difficult to handle. Numerically it would be possible
to have a decomposition $H\equiv C+D$ with $C\equiv
3L-\sum\limits_{i=1}^3x_i-
\sum\limits_{j=5}^{11}x_j-x_{12}$ which would be a plane cubic. In this case
$p(D)=4$, $H.D=6$. It might be interesting to check whether one can actually
construct surfaces with such a decomposition.

Of course, one can try and attempt to approach the problem of finding suitable
decompositions $H\equiv C+D$ more systematically. Let us assume $S$ is a
rational surface and $H\equiv C+D$ a decomposition to which the
Alexander-Bauer lemma can be applied. Let $h=h^1(\calO_S(H))$ be the
speciality of $S$. Since $C$ is mapped to a plane curve the exact sequence
$$
0\longrightarrow\calO_S(D)\longrightarrow\calO_S(H)\longrightarrow
\calO_C(H)\longrightarrow 0
$$
is exact on global sections, and hence
$$
h=h^1(D)+\delta(C)
$$
where $h^1(D)=h^1(\calO_S(D))$ and $\delta(C)=h^1(\calO_C(H))$. The analogous
sequence for $D$ and the assumption that $|H|$ restricts to a complete system
on the curves $D'\in|D|$ gives
$$
h=h^1(C)+\delta(D)
$$
where $h^1(C)$ and $\delta(D)$ are defined similarly. In general if $C$ is a
curve of genus $(d-1)(d-2)/2$ and $\calO_C(H)$ is a line bundle of degree $d$
it is difficult to show that $(C,\calO_C(H))$ is a plane curve. Hence it is
natural to assume $H.C\leq4$. In order to be able to control the linear system
$|H|$ on the curves $D'\in|D|$ one is normally forced to assume that $H.D\geq
2p(D)-2$ and $H|_D=K_D$ in case of equality. Hence $\delta(D)=0$ if
$H.D>2p(D)-2$ and $\delta(D)=1$ otherwise. Since $|H|$ is complete on $D$ we
have $h^0(\calO_D(H))\leq4$. Now using our assumption that $H.D\geq 2p(D)-2$
and Riemann-Roch on $D$ we find
$$
2p(D)-2\leq H.D\leq p(D)+3+\delta(D)
$$
and from this
$$
p(D)\leq 5+\delta(D).
$$
If $\delta(D)=0$ then $p(D)\leq5$. If $\delta(D)=1$ then $H|_D=K_D$ and
$h^0(\calO_D(H))=p(D)$, i.e.~$p(D)\leq4$ in this case. But now
$$
d=H.C+H.D\leq p(D)+7+\delta(D).
$$
This shows that one can find such a decomposition only if the degree $d\leq12$.
The case $d=12$ can only occur for $H.C=4$.

Finally we want to discuss the case $d=11$. In his thesis Popescu \cite{P}
gave three examples of rational surfaces of degree 11. In each case it is
$\PP^2$ blown up in 20 points. The linear systems are as follows:
\begin{eqnarray}
H&\equiv&10L-4x_1-3\sum_{i=2}^4x_i-2\sum_{j=5}^{14}x_j-\sum_{k=15}^{20}x_k
\label{IVgl1}\\
H&\equiv&11L-5x_1-3\sum_{i=2}^7x_i-2\sum_{j=8}^{13}x_j-\sum_{k=14}^{20}x_k
\label{IVgl2}\\
H&\equiv&13L-5x_1-4\sum_{i=2}^8x_i-2\sum_{j=9}^{11}x_j-\sum_{k=12}^{20}x_k
\label{IVgl3}
\end{eqnarray}
In each of these cases $S$ contains a plane quintic. The residual intersection
gives a pencil of rational (cases (\ref{IVgl1}) and (\ref{IVgl2})),
resp.~elliptic (case (\ref{IVgl3})) sextics. Since the linear system $|H|$ is
not complete on the curves of this linear system, one cannot immediately apply
the Alexander-Bauer lemma to this decomposition. One can ask whether there are
decompositions fulfilling the conditions given above. A candidate in case
(\ref{IVgl1}) is given by $C\equiv 4L-x_1-\sum\limits_{i=2}^4x_i-
\sum\limits_{j=5}^{14}x_j-\sum\limits_{k=15}^{17}x_k$ and $D\equiv H-C$. We do
not know whether surfaces with such a decomposition actually occur. In the
other cases one can show that no such decompositions exist.

\bibliographystyle{alpha}

\vspace{1cm}

Authors' addresses:

\bigskip

\parbox{5.5cm}{F.~Catanese\\ Dipartimento di Matematica\\ Universit\'a di
Pisa\\ via Buonarroti, 2\\ I 56127 Pisa\\ Italy\\
catanese\symbol{64}dm.unipi.it\\
catanese\symbol{64}vm.cnuce.cnr.it}
\parbox{5cm}{K.~Hulek\\ Institut f\"ur Mathematik\\ Universit\"at
Hannover\\ Postfach 6009\\ D 30060 Hannover\\ Germany\\
hulek$@$math.uni-hannover.de}
\end{document}